\documentclass{article}
\usepackage{amssymb}

\input{tcilatex}

\begin{document}

\begin{center}
\textbf{A\ SURVEY\ ON\ RECENT EXTENSIONS\ OF THE STIRLING FORMULA}%
\[
\]

Sorinel Dumitrescu$^{1}$ and Cristinel Mortici$^{2}$

\[
\]

$^{1}$Ph. D. Student, University Politehnica of Bucharest, Splaiul Independen%
\c{t}ei 313, Bucharest, Romania, sorineldumitrescu@yahoo.com\bigskip

$^{2}$Prof. dr. habil., Valahia University of T\^{a}rgovi\c{s}te, Bd. Unirii
18, 130082 T\^{a}rgovi\c{s}te, Romania, cristinelmortici@yahoo.com%
\[
\]
\end{center}

\[
\]

\begin{quote}
\textbf{ABSTRACT: }We present a survey on recent results about Stirling's
formula. More exactly, we reffer to a method using a form of Cesaro-Stolz
lemma firstly introduced in [C. Mortici Product approximations via
asymptotic integration Amer. Math. Monthly 117 (5) (2010) 434-441]. As an
example we improve a result obtained in [C. Mortici A substantial
improvement of the Stirling formula Appl. Math. Lett. 24 (2011) no. 8
1351-1354]. Finally, some numerical computations are made.
\end{quote}

\[
\]

\bigskip 

\textbf{MSC:} 2010: 26D15, 11Y25, 41A25, 34E05%
\[
\]

\textbf{INTRODUCTION}\bigskip

The factorial function$\ \ n!=1\cdot 2\cdot 3\cdot ...\cdot n$ \ (defined
for positive integers $n)$, and its extension gamma function%
\[
\Gamma (z)=\int_{0}^{\infty }t^{z-1}e^{-t}dt 
\]%
(to the real and complex values $z,$ excepting $-1,$ $-2,$ $-3,$ ...) has a
great importance in pure mathematics, as in applied mathematics and other
branches of science.

The Stirling formula [1], [14]%
\[
n!\sim \sqrt{2\pi n}\left( \frac{n}{e}\right) ^{n}:=\sigma _{n} 
\]%
was discovered by James Stirling (1692-1770) and Abraham de Moivre
(1667-1754) and represents one of the most used formulas in approximating
large factorials.

Although this formula is satisfactory in other branches such as engineering,
statistics, economics, in pure mathematics more accurate formulas are
required. In consequence, in the recent past, many authors gave increasingly
accurate formulas, but a sacrifice of simplicity. As such formulas are
implemented in some computer programs, they should be of a simple form in
order to avoid error cumulation.

Such simple formulas were obtained by authors starting from Stirling
formula. The following estimate slightly better than Stirling formula,%
\[
n!\sim \sqrt{2\pi }\left( \frac{n+\frac{1}{2}}{e}\right) ^{n+\frac{1}{2}}, 
\]%
was discovered by W. Burnside, see [4].

R. W. Gosper [7] introduced a much better approximation with a simple form%
\[
n!\sim \sqrt{2\pi \left( n+\frac{1}{6}\right) }\left( \frac{n}{e}\right)
^{n}=\gamma _{n}. 
\]

Similar results with Gosper formula were obtained by Mortici [13]%
\[
\sqrt{2\pi e}\cdot e^{-\omega }\left( \frac{n+\omega }{e}\right) ^{n+\frac{1%
}{2}}<n!<\sqrt{2\pi e}\cdot e^{-\zeta }\left( \frac{n+\zeta }{e}\right) ^{n+%
\frac{1}{2}}, 
\]%
where $\omega =\frac{3-\sqrt{3}}{6}$ and $\zeta =\frac{3+\sqrt{3}}{6}$. In
fact, $\omega $ and $\zeta $ are the values of $p$ which provide the best
possible approximations of the form%
\[
n!\sim \sqrt{2\pi e}\cdot e^{-p}\left( \frac{n+p}{e}\right) ^{n+\frac{1}{2}%
}. 
\]

Mortici [12] considered the family of approximations%
\begin{equation}
n!\sim \sqrt{2\pi n}\left( \frac{n}{e}+\frac{a}{n}\right) ^{n},\ \ \ a\in 
\mathbb{R}
\tag{1}
\end{equation}%
and proved that the most accurate result is obtained for $a=\frac{1}{12e}.$

As an illustration of this method, we propose in the next section the
following improvement of (1):%
\[
n!\sim \sqrt{2\pi \left( n+\frac{239}{181440n^{4}}\right) }\left( \frac{n}{e}%
+\frac{1}{12en}+\frac{1}{1440en^{3}}\right) ^{n}. 
\]%
\[
\]

\textbf{An example}\bigskip

We show in this section how the method using Cesaro-Stolz lemma works. As an
example, in the first part, we use the idea from [12] to improve Mortici
formula (1) to%
\[
n!\sim \sqrt{2\pi n}\left( \frac{n}{e}+\frac{1}{12en}+\frac{1}{1440en^{3}}%
\right) ^{n}:=w_{n}. 
\]%
In this sense, let us introduce the family of approximations%
\begin{equation}
n!\sim \sqrt{2\pi n}\left( \frac{n}{e}+\frac{a}{n}+\frac{b}{n^{3}}\right)
^{n}:=w_{n}\left( a,b\right) ,  \tag{2}
\end{equation}%
depending on real parameters $a,$ $b.$

The problem we rise here is what are the best parameters $a$ and $b$ which
provide the most accurate approximation $n!\sim w_{n}\left( a,b\right) .$

One method to measure the accuracy of an approximation formula of the form
(2) is to use the following\bigskip

\textbf{Lemma 2.1. }\emph{If }$\left( x_{n}\right) _{n\geq 1}$ \emph{is
convergent to zero} \emph{and}%
\[
\lim_{n\rightarrow \infty }n^{k}(x_{n}-x_{n+1})=l\in \lbrack -\infty ,\infty
], 
\]%
\emph{with }$k>1,$\emph{\ then}%
\[
\lim_{n\rightarrow \infty }n^{k-1}x_{n}=\frac{l}{k-1}. 
\]%
This lemma is a powerfull tool for accelerating some convergences, or for
constructing asymptotic series. See for proof and other details [2]-[13].

One way to measure the accuracy of an approximation of type (2) is to define
the sequence $\left( z_{n}\right) _{n\geq 1}$ by the relations%
\[
n!=\sqrt{2\pi n}\left( \frac{n}{e}+\frac{a}{n}+\frac{b}{n^{3}}\right)
^{n}\exp z_{n}\ ,\ \ \ n\geq 1 
\]%
and to consider an approximation of type (2) to be better as the sequence $%
\left( z_{n}\right) _{n\geq 1}$ converges faster to zero. In fact, we have%
\[
z_{n}=\ln n!-\ln \sqrt{2\pi }-\frac{1}{2}\ln n-n\ln \left( \frac{n}{e}+\frac{%
a}{n}+\frac{b}{n^{3}}\right) . 
\]%
In order to compute the speed of convergence of the sequence $z_{n}$ using
Lemma 2.1, we consider the difference%
\begin{eqnarray*}
z_{n}-z_{n+1} &=&-\ln \left( n+1\right) -\frac{1}{2}\ln \frac{n}{n+1} \\
&&-n\ln \left( \frac{n}{e}+\frac{a}{n}+\frac{b}{n^{3}}\right) +\left(
n+1\right) \ln \left( \frac{n+1}{e}+\frac{a}{n+1}+\frac{b}{\left( n+1\right)
^{3}}\right) .
\end{eqnarray*}%
By using a computer software such as Maple, we write $z_{n}-z_{n+1}$ as a
power series in $n^{-1},$%
\begin{eqnarray}
z_{n}-z_{n+1} &=&\frac{1}{n^{2}}\left( -ae+\frac{1}{12}\right) +\frac{1}{%
n^{3}}\left( ae-\frac{1}{12}\right) +\frac{1}{n^{4}}\left( -ae-3eb+\frac{3}{2%
}a^{2}e^{2}+\frac{3}{40}\right)  \TCItag{3} \\
&&+\frac{1}{n^{5}}\left( ae+6be-3a^{2}e^{2}-\frac{1}{15}\right)  \nonumber \\
&&+\frac{1}{n^{6}}\left( -ae-10be+5abe^{2}+5a^{2}e^{2}-\frac{5}{3}a^{3}e^{3}+%
\frac{5}{84}\right)  \nonumber \\
&&+O\left( \frac{1}{n^{7}}\right) \text{.}  \nonumber
\end{eqnarray}%
Now we are in a position to give the answer to the problem posed above. More
precisely, we formulate the following\bigskip

\textbf{Theorem 2.1. }\emph{(i) If }$a\neq \frac{1}{12e},$\emph{\ then the
rate of convergence of the sequence }$(z_{n})_{n\geq 1}$\emph{\ is }$n^{-1}$%
\emph{\ , since:}%
\[
\lim_{n\rightarrow \infty }nz_{n}=-ae+\frac{1}{12}\neq 0. 
\]%
\emph{(ii) If }$a=\frac{1}{12e},$\emph{\ and }$b\neq \frac{1}{1440e}$\emph{\
then the rate of convergence of the sequence }$(z_{n})_{n\geq 1}$\emph{\ is }%
$n^{-3}$\emph{\ , since:}%
\[
\lim_{n\rightarrow \infty }n^{3}z_{n}=\frac{1}{1440}-be\neq 0. 
\]%
\emph{(iii) If }$a=\frac{1}{12e},$\emph{\ and }$b=\frac{1}{1440e}$\emph{\
then the rate of convergence of sequence }$(z_{n})_{n\geq 1}$\emph{\ is }$%
n^{-5}$\emph{, since: }%
\[
\lim_{n\rightarrow \infty }n^{5}z_{n}=\frac{239}{362880}. 
\]%
As we explained before, the best approximation (2) is obtained in case
(iii), when the sequence $z_{n}$ is fastest possible, that is%
\[
n!\sim \sqrt{2\pi n}\left( \frac{n}{e}+\frac{1}{12en}+\frac{1}{1440en^{3}}%
\right) ^{n}:=w_{n}. 
\]

\emph{Proof of Theorem 2.1 }follows by Lemma 2.1.

(i) From (3) we get%
\[
\lim_{n\rightarrow \infty }n^{2}\left( z_{n}-z_{n+1}\right) =-ae+\frac{1}{12}
\]%
and by Lemma 2.1,%
\[
\lim_{n\rightarrow \infty }nz_{n}=-ae+\frac{1}{12}. 
\]%
(ii) and (iii). If $a=\frac{1}{12e},$ relation (3) becomes%
\[
z_{n}-z_{n+1}=\left( \frac{1}{480}-3be\right) \frac{1}{n^{4}}+\left( 6be-%
\frac{1}{240}\right) \frac{1}{n^{5}}+\left( \frac{361}{36\,288}-\frac{115}{12%
}be\right) \frac{1}{n^{6}}+O\left( \frac{1}{n^{7}}\right) . 
\]%
We have%
\[
\lim_{n\rightarrow \infty }n^{4}\left( w_{n}-w_{n+1}\right) =\frac{1}{480}%
-3be 
\]%
and%
\[
\lim_{n\rightarrow \infty }n^{3}w_{n}=\frac{1}{3}\left( \frac{1}{480}%
-3be\right) =\frac{1}{1440}-be. 
\]%
Finally, with $b=\frac{1}{1440e},$ we have%
\[
\lim_{n\rightarrow \infty }n^{6}\left( w_{n}-w_{n+1}\right) =\frac{239}{72576%
},\text{ \ \ and }\lim_{n\rightarrow \infty }n^{5}w_{n}=\frac{239}{362880}%
.\square 
\]%
Such improvements can continue similarly. Another idea is to introduce a new
real parameter $b$ and to consider the family of approximations%
\[
n!\sim \sqrt{2\pi \left( n+\frac{b}{n^{4}}\right) }\left( \frac{n}{e}+\frac{1%
}{12en}+\frac{1}{1440en^{3}}\right) ^{n},\text{ \ \ }b\in 
\mathbb{R}
\]%
together with the corresponding error sequence $t_{n}$ defined by%
\[
n!\sim \sqrt{2\pi \left( n+\frac{b}{n^{4}}\right) }\left( \frac{n}{e}+\frac{1%
}{12en}+\frac{1}{1440en^{3}}\right) ^{n}\exp t_{n}\ ,\ \ \ n\geq 1. 
\]%
As%
\[
t_{n}=\ln n!-\ln \sqrt{2\pi }-\frac{1}{2}\ln \left( n+\frac{b}{n^{4}}\right)
-n\ln \left( \frac{n}{e}+\frac{1}{12en}+\frac{1}{1440en^{3}}\right) , 
\]%
we use again Maple software to get%
\begin{eqnarray}
t_{n}-t_{n+1} &=&\left( -\frac{5}{2}b+\frac{239}{72\,576}\right) \frac{1}{%
n^{6}}+\left( \frac{15}{2}b-\frac{239}{24\,192}\right) \frac{1}{n^{7}} 
\TCItag{4} \\
&&+\left( -\frac{35}{2}b+\frac{26\,179}{1382\,400}\right) \frac{1}{n^{8}}%
+O\left( \frac{1}{n^{9}}\right) .  \nonumber
\end{eqnarray}%
The fastest sequence $t_{n}-t_{n+1}$ and consequently the fastest $t_{n}$
(see Lemma 2.1) are obtained when the first coefficient in (4) vanishes,
that is $b=\frac{239}{181\,440}.$ We deduce%
\begin{equation}
n!\sim \sqrt{2\pi \left( n+\frac{239}{181440n^{4}}\right) }\left( \frac{n}{e}%
+\frac{1}{12en}+\frac{1}{1440en^{3}}\right) ^{n}.  \tag{5}
\end{equation}%
Finally, we offer some computations which prove the superiority of formula
(5) over Mortici formula (1). Moreover, (5) is more accurate than Ramanujan
formula%
\[
n!\sim \sqrt{\pi }\left( \frac{n}{e}\right) ^{n}\sqrt[6]{8n^{3}+4n^{2}+n+%
\frac{1}{30}}. 
\]%
The next table contains the relative errors%
\[
\mu _{n}=\frac{\Gamma \left( n+1\right) }{\sqrt{2\pi n}\left( \frac{n}{e}+%
\frac{1}{12en}\right) ^{n}}-1 
\]%
\[
\rho _{n}=\frac{\Gamma \left( n+1\right) }{\sqrt{\pi }\left( \frac{n}{e}%
\right) ^{n}\sqrt[6]{8n^{3}+4n^{2}+n+\frac{1}{30}}}-1 
\]%
\[
\tau _{n}=\frac{\Gamma \left( n+1\right) }{\sqrt{2\pi \left( n+\frac{239}{%
181440n^{4}}\right) }\left( \frac{n}{e}+\frac{1}{12en}+\frac{1}{1440en^{3}}%
\right) ^{n}}-1 
\]%
\[
\begin{tabular}{|l|l|l|l|}
\hline
$n$ & $\mu _{n}$ & $\rho _{n}$ & $\tau _{n}$ \\ \hline
$10$ & $7.\,\allowbreak 003\,9\times 10^{-7}$ & $-8.\,\allowbreak
587\,2\times 10^{-8}$ & $-5.\,\allowbreak 795\,4\times 10^{-11}$ \\ \hline
$50$ & $5.\,\allowbreak 557\,5\times 10^{-9}$ & $-1.\,\allowbreak
496\,8\times 10^{-10}$ & $-7.\,\allowbreak 519\,1\times 10^{-16}$ \\ \hline
$100$ & $6.\,\allowbreak 945\times 10^{-10}$ & $-9.\,\allowbreak
451\,9\times 10^{-12}$ & $-5.\,\allowbreak 876\,8\times 10^{-18}$ \\ \hline
$500$ & $5.\,\allowbreak 555\,6\times 10^{-12}$ & $-1.\,\allowbreak
524\,7\times 10^{-14}$ & $-7.\,\allowbreak 524\,0\times 10^{-23}$ \\ \hline
\end{tabular}%
\]

\[
\]

\textbf{Acknowledgements. }The work of the second author was supported by a
grant of the Romanian National Authority for Scientific Research,
CNCS-UEFISCDI project number PN-II-ID-PCE-2011-3-0087.

\bigskip

\textbf{REFERENCES}\bigskip

[1]\ M. Abramowitz and I. A. Stegun, eds. Handbook of Mathematical Functions
with Formulas, Graphs, and Mathematical Tables, 9th printing, in : National
Bureau of Standards, Applied Mathematical Series, vol. 55, Dover, New York,
1972.

[2] V. Berinde and C. Mortici, New sharp estimates of the generalized
Euler-Mascheroni constant, Math. Inequal. Appl., 16 (2013), no. 1, 279-288.

[3] V. Berinde, A new generalization of Euler's constant, Creat. Math.
Inform. 18, 2 (2009), 123--128.

[4] W. Burnside, A rapidly convergent series for log\textit{N}!, Messenger
Math. 46(1917) 157-159

[5] Ch.-P. Chen, \textit{Inequalities and monotonicity properties for some
special functions}, J. Math. Inequal. \textbf{3} (2009), 79--91.

[6] Ch.-P. Chen, \textit{The Best Bounds in Vernescu's Inequalities for the
Euler's Constant}, RGMIA Res. Rep. Coll. \textbf{12} (2009), no.3, Article
11.[21] C. Mortici, An ultimate extremely accurate approximation formula for
the factorial function, Arch. Math. (Basel) 93(1)(2009) 37-45.

[7] R. W. Gosper, Decision procedure for indefinite hypergeometric
summation, Proc. Natl. Acad. Sci. USA 75(1918) 40-42

[8] Hu Yue, A strenghtened Carleman's inequality, Commun. Math. Anal., 1
(2006), no. 2, 115-119.

[9] C. Mortici, New approximations of the gama function in terms of the
digamma function, Appl. Math. Lett. 23 (1) (2010) 97-100.

[10] C. Mortici, Optimizing the rate of convergence in some new classes of
sequences convergent to Euler's constant, Anal. Appl. (Singap.) 8 (1) (2010)
99-107

[11] C. Mortici, Product approximations via asymptotic integration, Amer.
Math. Monthly 117 (5) (2010) 434-441.

[12] C. Mortici, A substantial improvement of the Stirling formula, Appl.
Math. Lett., 24 (2011), no. 8, 1351-1354.

[13] C. Mortici, An ultimate extremely accurate approximation formula for
the factorial function, Arch. Math. (Basel) 93(1)(2009) 37-45.

[14] J. Stirling, Methodus Differentialis, Sive Tractatus de Summation of
Interpolation Serierum Infinitarium, London, 1730. English translation by J.
Holliday, The Differential Method: A Treatise of the Summation and
Interpolation of Infinite Series. \ James Stirling's Methodus
Differentialis: An Annotated Translation of Stirling's text, Sources and
Studies in the History of Mathematics and Physical Science, Springer -
Verlag, London, 2003.

\end{document}